\documentclass[12pt]{amsart}
\usepackage{amscd}

\theoremstyle{definition}

\def \ph{\varphi}

\def \refeq#1{equation (\ref{#1})}
\def \ra{\rightarrow}

\def \hom{\mbox{\rm Hom}}

\def \tns{\otimes}

\def \mplus{+\cdots+}
\def \mcom{,\cdots,}
\def \k{\mbox{$\mathbb K$}}

\def \C{\mbox{$\mathbb C$}}
\def \Z{\mbox{$\mathbb Z$}}

\def\br#1#2{\lbrack#1,#2\rbrack}

\def\zt{\mbox{$\Z_2$}}

\def\ad{\operatorname{ad}}

\def\d{d}
\def\td{\tilde\d}

\def\A{\mbox{$\mathcal A$}}
\def\B{\mbox{$\mathcal B$}}
\def\L{L}

\def\LA{\mbox{$\L_{\A}$}}
\def\LB{\mbox{$\L_{\B}$}}

\def\m{\mbox{$\mathfrak m$}}

\def\coder{\operatorname{Coder}}

\def\linf{\mbox{$L_\infty$}}
\def\and{\mbox{ \rm and }}

\def\s#1{(-1)^{#1}}
\DeclareMathOperator*{\invlim}{\overleftarrow{\rm lim}}

\def\htns{\mbox{$\hat\tns\,\,$}}

\def\phd#1#2{\ph^{#1}_{#2}}

\def\psd#1#2{\psi^{#1}_{#2}}

\def\parity{\operatorname{parity}}

\author{Alice Fialowski}
\address{E\"otv\"os Lor\'and University\\
Budapest, Hungary} \email{fialowsk@cs.elte.hu}
\author{Michael Penkava}
\address{University of Wisconsin\\
Eau Claire, WI 54702-4004}
\email{penkavmr@uwec.edu}
\subjclass{14D15,13D10,14B12,16S80,16E40,\\17B55,17B70}
\keywords{Versal Deformations, Infinity Algebras, Cohomology,
Infinitesimal Deformation, Extensions}
\thanks{The research of the authors was supported by grants
MTA-OTKA-NSF 38453, OTKA T043641 and T043034 and by grants
from the University of Wisconsin-Eau Claire} \title[Examples of
Miniversal Deformations of Infinity Algebras]{Examples of Miniversal
Deformations of Infinity Algebras}

\begin{document}
\setlength{\multlinegap}{0pt}
\begin{abstract}
A classical problem in algebraic deformation theory is whether an
infinitesimal deformation can be extended to a formal deformation. The
answer to this question is usually given in terms of Massey powers.
If all Massey powers of the cohomology class determined by the
infinitesimal deformation vanish, then the deformation extends to a
formal one.  We consider another approach to this problem, by
constructing a miniversal deformation of the algebra. One advantage of
this approach is that it answers not only the question of existence,
but gives a construction of an extension as well.
\end{abstract}
\date\today
\maketitle

In this paper, we study some examples of miniversal deformations of
infinity algebras, and use these examples to illustrate how to use a
miniversal deformation to determine when an infinitesimal deformation
extends to a formal deformation. Actually, using a miniversal
deformation, one can construct such an extension explicitly. Also, the
obstruction to an extension can be computed by this method.

An infinitesimal deformation extends to a formal one precisely when
the unique morphism from the base of the universal infinitesimal
deformation to the base of the given deformation inducing the
infinitesimal deformation can be lifted to a morphism from the
miniversal deformation to the formal power series ring in the
parameter of the deformation.  A nice property of this algebraic
approach is that it answers more than the question of existence; in
fact, it gives a construction of an extension of the infinitesimal
deformation to a formal deformation. Moreover, since the question is
reduced to studying the morphisms of the base of the miniversal
deformation to a formal power series ring, the general question of
which infinitesimal deformations extend to formal deformations is
reduced to a simple algebraic question.

Let us point out that the problem of extending a deformation emerges in
deformation quantization as well, see \cite{cft}.

\section{Introduction} We work in the framework of the parity
reversion $W=\Pi V$ of the usual vector space $V$ on which an
\linf\ algebra structure is defined,  because in the $W$
framework,  an \linf\ structure is simply an odd coderivation $d$
of the symmetric coalgebra $S(W)$, satisfying $d^2=0$,  in other
words, it is an odd codifferential in the \zt-graded Lie algebra
of coderivations of $S(W)$.  As a consequence, when studying
\zt-graded Lie algebra structures on $V$, the parity is reversed,
so that an $m|n$-dimensional vector space $W$ corresponds to a
$n|m$-dimensional \zt-graded Lie structure on $V$.  Moreover,  the
\zt-graded \emph{anti-symmetry} of the Lie bracket on $V$ becomes
the \zt-graded \emph{symmetry} of the associated coderivation $d$
on $S(W)$.

A formal power series $d=d_1+\cdots$, with $d_i\in L_i=\hom(S^i(W),W)$
determines an element in $L=\hom(S(W),W)$, which is naturally
identified with $\coder(S(W))$, the space of coderivations of the
symmetric coalgebra $S(W)$. Thus $L$ is a \zt-graded Lie algebra. An
odd element $d$ in $L$ is a called a  \emph{codifferential} if $\br
dd=0$. We also say that $d$ is an \linf\ structure on $W$.

A  detailed description of \linf\ algebras can be obtained in
\cite{pen3,pen4}. The study of examples of \linf\ algebra structures
in \cite{fp2,fp3,fp4},and especially \cite{bfp1} may be useful to the
reader because they contain many examples of \linf\ algebras and their
miniversal deformations.

Let us establish some basic notation for the cochains. Suppose
$W=\langle w_1\mcom w_{m+n}\rangle$ with $w_1\mcom w_n$ odd and
$w_{n+1}\mcom w_{m+n}$ even elements. If $I=\{i_1\mcom i_{m+n}\}$ is a
multi-index, with $i_k$ either zero or one when $k\le n$, let
$w_I=w_1^{i_1}\cdots w_{m+n}^{i_{m+n}}$.
Denote $\deg(I)=i_1\mplus i_{m+n}$, $\parity(I)= i_1\mplus i_n\pmod
2$. Then for $n\ge 1$,
\begin{align*}
(S^{n}(W))_e=&\langle w_I |\parity(I)=0, \deg(I)=n\rangle\\
(S^{n}(W))_o=&\langle w_I |\parity(I)=1, \deg(I)=n\rangle
\end{align*}
For $n=\deg(I)$, for  $j=1..m+n$, we  define a map
$\ph^I_j:S^n(W)\ra W$ by $\ph^I_j(w_J)=I!\delta^I_J w_j$, where
$I!=i_1!\cdots i_{m+n}!$. Let $L_n:=\hom(S^n(W),W)$, then
$L_n=\langle \ph^I_j, \deg(I)=n \rangle$. If $\ph$ is odd, we
denote it by the symbol $\psi$ to make it easier to distinguish
the even and odd elements.

\subsection{ Versal Deformations}
For a treatment of classical formal deformation theory we refer to
\cite{gers1}. Versal deformation theory was first worked out for
the case of Lie algebras in \cite{fi1, fi2, ff2} and then extended
to \linf\ algebras in \cite{fp1}.

An augmented local ring $\A$ with maximal ideal $\m$ will be
called an \emph{infinitesimal base} if  $\m^2=0$, and a
\emph{formal base} if $\A=\invlim_n \A/\m^n$. A deformation of an
\linf\ algebra structure $d$ on $W$ with  base given by a local
ring $\A$ with augmentation $\epsilon:\A\ra\k$, where $\k$ is the
field over which $W$ is defined, is an $\A$-\linf\ structure $\td$
on $W\htns\A$ such that the morphism of $\A$-\linf\ algebras
$$\epsilon_*=1\tns\epsilon:\LA=L\tns \A\ra L\tns\k=L$$ satisfies
$\epsilon_*(\td)=d$.  (Here $W\htns\A$ is an appropriate
completion of $W\tns\A$.) The deformation is called \emph{infinitesimal}
(\emph{formal}) if $\A$ is an infinitesimal (formal) base.

In general, the cohomology $H(D)$ of $d$ given by the operator $D:L\ra
L$ with $D(\ph)=\br\ph d$ may not be finite dimensional. However, $L$
has a natural filtration $L^n=\prod_{i=n}^\infty L_i$, which induces a
filtration $H^n$ on the cohomology, because $D$ respects the
filtration. We say that $H(D)$ is of \emph{finite type} if $H^n/H^{n+1}$ is finite
dimensional for all $n$. Since this is always true when $W$ is finite dimensional,
the examples we study here will always be of finite type.  A set
$\{\delta_i\}$ will be called a \emph{basis of the cohomology}, if any
element $\delta$ of the cohomology can be expressed uniquely as a
formal sum $\delta=\delta_i a^i$. (Here and throughout the paper, we use Einstein's summation convention).
If we identify $H(D)$ with a
subspace of the space of cocycles $Z(D)$, and we choose a \emph{basis}
$\{\beta_i\}$ of the coboundary space $B(D)$, then any element
$\zeta\in Z(D)$ can be expressed uniquely as a sum $\zeta=\delta_i a^i
+\beta_i b^i$.

For each $\delta_i$, let $u^i$ be a parameter of opposite parity.
Then the infinitesimal deformation $d^1=d+\delta_i u^i$, with base
$\A=\k[u^k]/(u^iu^j)$ is universal in the sense that if $\d^i$ is any
infinitesimal deformation with base $\B$, then there is a unique
morphism $f:\A\ra\B$, such that the morphism $f_*=1\tns
f:\LA\ra\LB$ satisfies $f_*(d^1)\sim d^i$.

For formal deformations, there is no universal object in the sense
above. A \emph{versal deformation} is a deformation $d^\infty$ with
formal base $\A$ such that if $d^f$ is any formal deformation with
base $\B$, then there is some morphism $f:\A\ra\B$ such that
$f_*(d^\infty)\sim d^f$. If $f$ is unique whenever $\B$ is
infinitesimal, then the versal deformation is called
\emph{miniversal}. In \cite{fp1}, we constructed a miniversal
deformation for \linf\ algebras with finite type cohomology.

The method of construction is as follows.  Define a coboundary
operator $D$ by $D(\ph)=[\ph,d]$. First, one constructs the universal
infinitesimal deformation $d^1=d+\delta_i u^i$, where $\delta_i$ is a
graded basis of the cohomology $H(D)$ of $d$, or more correctly, a
basis of a subspace of the cocycles which projects isomorphically to a
basis in cohomology, and $u^i$ is a parameter whose parity is opposite
to $\delta_i$.  The infinitesimal assumption that the products of
parameters are equal to zero gives the property that $[d^1,d^1]=0$.
Actually, we can express
\begin{equation*}
[d^1,d^1]=\s{\delta_j(\delta_i+1)}[\delta_i,\delta_j]u^iu^j=\delta_k
a^k_{ij}u^iu^j +\beta_k b^k_{ij}u^iu^j,
\end{equation*}
where $\beta_i$ is a basis of the coboundaries, because the
bracket of $d^1$ with itself is a cocycle. Note that the right
hand side is of degree 2 in the parameters, so it is zero up to
order 1 in the parameters.

If we suppose that $D(\gamma_i)=-\frac12\beta_i$, then by
replacing $d^1$ with $$d^2=d^1+\gamma_kb^k_{ij}u^iu^j,
$$ one obtains
\begin{equation*}
[d^2,d^2]=\delta_k a^k_{ij}u^iu^j+ 2[\delta_l
u^l,\gamma_kb^k_{ij}u^iu^j]+[\gamma_kb^k_{ij}u^iu^j,\gamma_lb^l_{ij}u^iu^j]
\end{equation*}
Thus we are able to get rid of terms of degree 2 in the
coboundary terms $\beta_i$, but those which involve the cohomology
terms $\delta_i$ can not be eliminated.  Therefore,
$R^k=a^k_{ij}u^iu^j$ must be equal to zero up to order 3, which is accomplished by taking the base
of the second order deformation to be the  quotient of the ring $\k[u^l]/(u^iu^ju^k)$ by the
ideal generated by the \emph{second order relations}
$R^k$.  One continues this
process, taking the bracket of the \emph{$n$-th order deformation}
$d^n$, adding some higher order terms to cancel coboundaries,
obtaining higher order relations, which extend the second order
relations.

Either the process continues indefinitely, in which
case the miniversal deformation is expressed as a formal power
series in the parameters, or after a finite number of steps, the
right hand side of the bracket is zero after applying the $n$-th
order relations. In this case, the miniversal deformation is
simply the $n$-th order deformation. In either case, we obtain a set
of relations $R^i$ on the parameters, one for each $\delta_i$, and
the algebra $A=\C[[u^i]]/(R^i)$ is called the base of the
miniversal deformation. Examples of the construction of miniversal
deformations can be found in \cite{ff3,ff2,fipo,fp2,fp3}.
\subsection{Extensions of Infinitesimal Deformations}
Let us put together a general picture of how to use a miniversal
deformation to solve the extension problem. Let us suppose that
$\psi_k$ and $\phi_k$ are bases of the odd and even parts of a
preimage of the cohomology of a codifferential $d$, and $\alpha_k$,
$\beta_k$ are bases of the odd and even parts of a preimage of the
coboundaries determined by $d$. Then there is a miniversal deformation
of the form
\begin{equation*}
d^\infty=d+\psi_k t^k +\phi_k \theta^k + \alpha_k x^k +\beta_k y^k,
\end{equation*}
where the $t^k$ are odd parameters, the $\theta^k$ are even ones,
the $x^k$ are odd and the $y^k$ are even formal power series in the
parameters, and for each $k$, there are odd relations $r_o^k$ and even
relations $r_e^k$, which are formal power series in the parameters.
The base of the miniversal deformation is given by
$\k[[t^k,\theta^l]]/(r_o^k,r_e^l)$.

Classically, an infinitesimal deformation is given by a single even
parameter $u$. It is natural to extend the classical picture by adding
an odd parameter $\theta$, so that for our purposes we will state the
deformation problem in the following manner. Consider an infinitesimal
deformation of the form
\begin{equation*}
d^i=d+\psi u +\ph \theta,
\end{equation*}
where $\psi$ is an odd and $\ph$ is an even cocycle.  An important
question is:

\emph{When does this infinitesimal deformation extend to a
formal deformation?}

Without loss of generality, one can assume that $\psi$ and $\ph$ are
(possibly infinite) linear combinations of the $\psi_k$ and $\phi_k$.
This is because one can remove any coboundary term by applying an
equivalence.  Similarly, any extension of this infinitesimal
deformation to an $n$-th order deformation is equivalent to one of the
form
\begin{equation*}
d^n=d+\psi_ka^k_iu^i+\phi_k b^k_iu^i\theta +\alpha_k g^k_iu^i +\beta_k h^k_iu^i\theta,
\end{equation*}
where $a^k=a^k_iu^i$, $b^k=b^k_iu^i$, $g^k=g^k_iu^i$ and
$h^k=h^k_iu^i$ are polynomials of degree less than or equal to $n$ in
$u$ without constant term such that $\psi=\psi_k a_1^k$ and
$\ph=\phi_kb_1^k$. An important generalization of the first question is:

\emph{When does such an
$n$-th order deformation extend to a formal deformation?}

To answer this question, first note that if we identify $t^k=a^k$ and
$\theta^k=b^k\theta$, then
\begin{enumerate}
\item The relations on the base are satisfied
up to order $n+2$.
\item $g^k=x^k \pmod{n+1}$ and $h^k\theta= y^k\pmod{n+1}$.
\end{enumerate}

The deformation $d^n$ extends to a formal deformation $d^f$ of $d$
precisely when there are extensions of $a^k$ and $b^k$ to formal power
series such that the identifications $t^k=a^k$ and $\theta^k=b^k$
satisfy the relations on the base of the formal deformation. Let $f$ be the morphism
$f:\k[[t^k,\theta^k]]\ra \k[[u,\theta]]$ induced by the
identifications above. Then $f$ descends to a morphism from the base of
the miniversal deformation to $\k[[u,\theta]]$, and
$d^f=f_*(d^\infty)$ is a formal deformation extending $d^n$. Because
there may be many extensions of $a^k$ and $b^k$ to formal power series
satisfying the relations, the deformation $d^f$ is not unique in
general.

Given a formal deformation $d^f$ of the form
\begin{equation}\label{formaldef}
d^f=d+\psi_ka^k_iu^i+\phi_k b^k_iu^i\theta +\alpha_k g^k_iu^i +\beta_k h^k_iu^i\theta,
\end{equation}
where now $a^k$, $b^k$, $g^k$ and $h^k$ are formal power series, there
is a unique map $f$ from the base of the miniversal deformation to
$\k[[u,\theta]]$ satisfying $f_*(d^\infty)=d^f$. Thus it may seem that
the miniversal deformation is universal.  The problem is that we work
in the category of equivalence classes of deformations, so that $d^f$
may be equivalent to other deformations of the form given by
\refeq{formaldef}. The necessity of working with equivalence classes
is clear from the fact that in general, a formal deformation is not of
the form given by \refeq{formaldef}, but merely equivalent to one in
such a form, because coboundary terms may appear in $d^f$. Moreover,
an equivalent deformation in the form given by \refeq{formaldef} is not unique
in general. We will
give an example later on in the text to illustrate this point.

Our purpose in this article is to construct some nontrivial examples
of miniversal deformations and use them to illustrate how to carry out
the procedure of determining which infinitesimal deformations extend
to a formal one.

\section{Codifferentials on a $2|1$ dimensional space}
In this section, we will be studying miniversal deformations of some
\linf\ structures on a $2|1$ dimensional space. For this space, our
multiindices $I$ will be ordered triples. Since $w_1$ is the only odd
basis element, we have
\begin{align*}
(L_{n})_e&=\langle\ph^{1,q,n-q-1}_1,\ph^{0,p,n-p}_2,\ph^{0,p,n-p}_3|1\le q\le n-1,1\le p\le n\rangle\\
(L_{n})_o&=\langle\psi^{1,q,n-q-1}_2,\psi^{1,q,n-q-1}_3,\psi^{0,p,n-p}_1|1\le q\le n-1,1\le p\le n\rangle,
\end{align*}
so that $|L_n|=3n+2|3n+1$. In \cite{bfp1}, the moduli space of
codifferentials of degree two on this space was computed.  The degree
two codifferentials are divided into two kinds.  Codifferentials of
the \emph{first kind} are of the form
\begin{equation*}
d=\psd{1,1,0}2x+\psd{1,1,0}3a+\psd{1,0,1}2b
+\psd{1,0,1}3c
\end{equation*}
and those of the \emph{second kind} are of the form
\begin{equation}\label{second}
d=\psd{0,2,0}1a+ \psd{0,1,1}1b + \psd{0,0,2}1c.
\end{equation}
Miniversal deformations for codifferentials of the first kind were
computed in \cite{bfp1}. The codifferentials of degree two of the
first kind form a complicated one parameter family, while the
codifferentials of degree two of the second kind are all equivalent to one of only two
types, which we called Type $(1,0,0)$ and Type $(0,1,0)$, where the
type represents the triple $(a,b,c)$ of coefficients in
\refeq{second}. Even though the description of the moduli space of
degree two codifferentials of the second kind is simple, the cohomology for both of the
degree two codifferentials of the second kind is infinite dimensional.
We did not give a complete description of the miniversal deformations
of degree two codifferentials of the second kind in \cite{bfp1}, so we
will give that description here.  The miniversal deformations provide
some nice examples which illustrate how to use a miniversal
deformation to determine whether an infinitesimal deformation extends
to a formal deformation.

\subsection{Miniversal deformations of Type (1,0,0)}
Let $d=\psd{0,2,0}1$.
We obtain the following table
of coboundaries.
\begin{align*}
D(\phd{1, q, n-q-1}1) &=\psd{0, 2+q, n-q-1}1\\
D(\phd{0, p, n-p}2) &=-2\psd{0,p+1, n-p}1\\
D(\phd{0, p, n-p}3) &=0\\
D(\psd{0, p, n-p}1) &=0\\
D(\psd{1, q, n-q-1}2) &=2\phd{1,q+1,n-q-1}1+\phd{0,q+2,n-q-1}2\\
D(\psd{1, q, n-q-1}3) &=\phd{0,q+2,n-q-1}3
\end{align*}
The cohomology is given by
\begin{align*}
&H^1=\langle \psd{0,0,1}1,\psd{0,1,0}1,\phd{0,0,1}3,\phd{0,1,0}3,2\phd{1,0,0}1+
    \phd{0,1,0}2\rangle\\
&H^n=\langle \psd{0,0,n}1,\phd{0,0,n}3,\phd{0,1,n-1}3,2\phd{1,0,n-1}1+
    \phd{0,1,n-1}2,\rangle,\quad\text{if $n>1$}
\end{align*}

Let us label the cohomology classes as follows
\begin{gather*}
\xi=\psd{0,1,0}1\\
\psi_n=\psd{0,0,n}1,\quad
\phi_n=\phd{0,0,n}3,n>0\\
\sigma_n=\phd{0,1,n-1}3,\quad\tau_n=2\phd{1,0,n-1}1+\phd{0,1,n-1}2,\quad n> 0
\end{gather*}
In order to construct the miniversal deformation,  we choose
pre-images of the coboundaries as follows:
\begin{equation*}
\gamma_{k,l}=\tfrac12\phd{0,k,l}2,\quad
\alpha_{k,l}=\psd{1,k,l}2,\quad
\beta_{k,l}=\psd{1,k,l}3
\end{equation*}
Then
\begin{align*}
D(\gamma_{k,l})=&-\psd{0,k+1,l}1\quad
D(\alpha_{k,l})=&2\phd{1,k+1,l}1+\phd{0,k+2,l}2\quad
D(\beta_{k,l})=&\phd{0,k+2,l}3
\end{align*}
The universal infinitesimal deformation is given by
\begin{equation*}
d^1=\psd{0,2,0}1+\xi s^1 +\psi_n t^n+\phi_n \theta^n
 +\sigma_n\eta^n+\tau_n\zeta^n,
\end{equation*}
where $s^1$ and $t^n$ are even parameters and $\theta^n$, $\eta^n$ and
$\zeta^n$ are odd parameters.

The brackets we need to compute $[d^1,d^1]$ are
\begin{align*}
[\xi,\phi_k]=&[\xi,\sigma_k]=[\tau_k,\tau_l]=0\\
[\xi,\tau_k]=&
D(\gamma_{0,k-1})
&[\phi_k,\phi_l]=&\phi_{k+l-1}(k-l)\\
[\psi_k,\phi_l]=&\psi_{k+l-1}k
&[\phi_k,\sigma_l]=&\sigma_{k+l-1}(k-l+1)\\
[\psi_k,\sigma_l]=&
-D(\gamma_{0,k+l-2})k
&[\phi_k,\tau_l]=&\tau_{k+l-1}(1-l)\\
[\psi_k,\tau_l]=&-2\psi_{k+l-1}
&[\sigma_k,\sigma_l]=&
D(\beta_{0,k+l-3})(k-l)\\
[\sigma_k,\tau_l]=&
\sigma_{k+l-1}+D(\alpha_{0,k+l-3})(1-l)
\end{align*}
Note that there are two exceptions to the rules above:
\begin{align*}
[\xi,\tau_1]=-\xi,\quad
[\psi_1,\sigma_1]=\xi
\end{align*}
The second order relations are given by
\begin{gather*}
-s^1\zeta^1+t^1\eta^1=0\\
\sum_{k+l=n+1} t^k (k\theta^l-2\zeta^l)=
\tfrac12\sum_{k+l=n+1}(k-l)\theta^k\theta^l=0\\
\sum_{k+l=n+1}\eta^k(\zeta^l+(k-l-1)\theta^l)=
\sum_{k+l=n+1}(1-l)\theta^k\zeta^l=0
\end{gather*}

The second order deformation is easily computed to be
\begin{equation*}
d^2=d^1+\gamma_{0,l}x^l+\alpha_{0,l}y^l+\beta_{0,l}z^l,
\end{equation*}
where
\begin{align*}
x^n=&s^1\zeta^{n+1}-\sum_{k+l=n+2}kt^k\eta^l\\
y^n=&\sum_{k+l=n+3}(l-1)\eta^k\zeta^l\\
z^n=&\tfrac12\sum_{k+l=n+3}(l-k)\eta^k\eta^l=-\sum_{k+l=n+3}k\eta^k\eta^l
\end{align*}
Note that only some of the pre-images of coboundaries actually play
any role in the second order deformation. In this example, it turns
out that the second order deformation is miniversal, so these are the
only cochains which are necessary to add. To see this, let us consider
the brackets which arise in the computation of $[d^2,d^2]$.
\begin{align*}
[\xi,\alpha_{0,n}]=&\tfrac12\tau_{n+1}+\gamma_{1,n}&
[\xi,\gamma_{0,n}]=&\tfrac12\psi_n\\
[\psi_k,\alpha_{0,l}]=&2\gamma_{0,k+l}&
[\psi_k,\gamma_{0,l}]=&0\\
[\alpha_{0,k},\phi_l]=&\alpha_{0,k+l-1}k&
[\phi_k,\gamma_{0,l}]=&-\gamma_{0,k+l-1}l\\
[\alpha_{0,k},\sigma_l]=&\alpha_{1,k+l-2}k-\beta_{0,k+l-1}&
[\sigma_k,\gamma_{0,l}]=&\tfrac12\phi_{k+l-1}-\gamma_{1,k+l-2}l\\
[\alpha_{0,k},\tau_l]=&\alpha_{0,k+l-1}&
[\tau_k,\gamma_{0,l}]=&\gamma_{0,k+l-1}\\
[\xi,\beta_{0,n}]=&\sigma_{n+1}&
[\alpha_{0,k},\alpha_{0,l}]=&0\\
[\psi_k,\beta_{0,l}]=&\phi_{k+l}+\tau_{k+l}\tfrac k2-\gamma_{1,k+l-1}k&
[\alpha_{0,k},\beta_{0,l}]=&0\\
[\beta_{0,k},\phi_l]=&\beta_{0,k+l-1}(k-l)&
[\beta_{0,k},\beta_{0,l}]=&0\\
[\beta_{0,k},\sigma_l]=&\beta_{1,k+l-2}(k+1-l)&
[\alpha_{0,k},\gamma_{0,l}]=&0\\
[\beta_{0,k},\tau_l]=&\alpha_{1,k+l-2}(1-l)+2\beta_{0,k+l-1}&
[\gamma_{0,k},\gamma_{0,l}]=&0\\
[\beta_{0,k},\gamma_{0,l}]=&-\alpha_{0,k+l-1}\left(\tfrac l2\right)
\end{align*}
No coboundaries appear in these brackets, which means that the second
order deformation is miniversal. It is also the case that the sum of
the terms in the bracket $[d^2,d^2]$ involving the pre-images of a coboundary must vanish. In
particular, the sum of all terms involving
$\alpha_{1,n}$, $\beta_{1,n}$ or $\gamma_{1,n}$ cochains must vanish.
Strictly speaking, it is unnecessary to check this fact, since it is
guaranteed by the existence theorem for the miniversal deformation
\cite{fp1}, but we found it interesting to check the manner in which
the terms cancel. In fact,  these terms cancel without using the
relations on the base, although, as we shall show later, the same is
not true for the $\alpha_{0,n}$, $\beta_{0,n}$ and $\gamma_{0,n}$
cochains.

The $\beta_{1,n}$ terms appear only in the bracket $[\beta_{0,k-1},
\sigma_l]z^{k-1}\eta^l$, so we should have
$\sum_{k+l=n+2}(k+1-l)z^{k-1}\eta^l=0$. We obtain
\begin{equation*}
\sum_{k+l=n+3}\!\!(k-l)z^{k-1}\eta^l=\!\!\sum_{i+j+l=n+5}\!\!-
    (i+j-2-l)i\eta^i\eta^j\eta^l,
\end{equation*}
which vanishes simply because the $\eta$ cochains anti-commute.

To see that the $\alpha_{1,n}$ cochains cancel, note that there are
two sources of such terms. From
$[\beta_{0,k-1},\tau_l]z^{k-1}\zeta^l$, we get
$\alpha_{1,k+l-3}(1-l)z^{k-1}\zeta^l$, while from
$[\alpha_{0,k-1},\sigma_l]y^{k-1}\eta^l$, we get
$\alpha_{1,k+l-3}ky^{k-1}\eta^l$. Substituting for $y^{k-1}$ and
$z^{k-1}$, and summing, we obtain
\begin{equation*}
\sum_{i+j+l=n+5}\!\!\!-(1-l)i\eta^i\eta^j\zeta^l-
    (i+j-2)(1-j)\eta^i\zeta^j\eta^l=\!\!\!
\sum_{i+j+l=n+5}\!\!l(2-l)\eta^i\eta^j\zeta^l=0
\end{equation*}

To see that the $\gamma_{1,n}$ terms vanish, we compute
\begin{align*}
[\sigma_k,\gamma_{0,l}]\eta^k x^{l}=&-\gamma_{1,k+l-2}\eta^k x^{l}\\
[\psi_k,\beta_{0,l-1}] t^kz^{l-1}=&-\gamma_{1,k+l-2}kt^kz^{l-1}\\
[\xi,\alpha_{0,n}]s^1y^n=&\gamma_{1,n}s^1y^n
\end{align*}
Therefore, the sum of all terms involving $\gamma_{1,n}$ has coefficient
\begin{align*}
s^1y^n-\sum_{k+l=n+2}l\eta^ks^1\zeta^{l+1`}+
    \sum_{k+i+j=n+4}(n+2-k)\eta^kit^i\eta^j+kt^ki\eta^i\eta^j=0,
\end{align*}
since the first sum above is just $s^1y^n$ and the second sum vanishes by the
anticommutativity of $\eta$ cochains.

The relations on the base of the miniversal deformation are
\begin{gather*}
r_1=-s^1\zeta^1+t^1\eta^1=0\\
r_2^n=\tfrac12s^1x^n+\sum_{k+l=n+1} t^k (k\theta^l-2\zeta^l)=0\\
r_3^n=\sum_{k+l=n+1}\tfrac12(k-l)\theta^k\theta^l+t^kz^{l-1}+
    \tfrac12\eta^kx^l=0\\
r_4^n=s^1z^{n-1}+\sum_{k+l=n+1}\eta^k(\zeta^l+(k-l-1)\theta^l)=0\\
r_5^n=\tfrac12s^1y^{n-1}+\sum_{k+l=n+1}(1-l)\theta^k\zeta^l+
    \tfrac k2t^kz^{l-1}=0
\end{gather*}
Now let us show that the coefficients of the $\alpha_{0,n}$,
$\beta_{0,n}$ and $\gamma_{0,n}$ cochains vanish.

The coefficients of the terms involving $\beta_{0,n}$ are
\begin{align*}
[\alpha_{0,k-1},\sigma_l]y^{k-1}\eta^l=&-\beta_{0,k+l-2}y^{k-1}\eta^l\\
[\beta_{0,k-1},\phi_l]z^{l-1}\theta^l=&\beta_{0,k+l-2}(k-1-l)z^{k-1}\theta^l\\
[\beta_{0,k-1},\tau_l]z^{l-1}\zeta^l=&\beta_{0,k+l-2}2z^{k-1}\zeta^l
\end{align*}
First, we observe that
\begin{equation*}
\sum_{k+l=n+2}y^{k-1}\eta^l=\sum_{i+j+l=n+4}(j-1)\eta^i\zeta^j\eta^l=0,
\end{equation*}
so the coefficients from the $[\alpha_{0,k-1},\sigma_l]$ terms add up
to zero on their own. Next, we have
\begin{equation*}
\sum_{i+j=n+2}(j-i)\eta^iz^{j-1}=-\sum_{i+k+l=n+4}(n+2-2i)\eta^ik\eta^k\eta^l=0.
\end{equation*}
Thus we have
\begin{align*}
0=\sum_{i+k=n+3}(k-i)\eta^ir_4^j=&
\sum_{i+k+l=n+4}(j+l-1-i)\eta^i\eta^j(\zeta^l+(j-l-1)\theta^l)\\=&
\sum_{i+k+l=n+4}-2i\eta^i\eta^j\zeta^l+(-i^2+(l+3)i)\eta^i\eta^j\theta^l\\=&
-\sum_{i+j+l=n+4}i\eta^i\eta^j((i+j-3-l)\theta^l+2\zeta^l)\\=&
\sum_{k+l=n+2}z^{k-1}((k-1-l)\theta^l+2\zeta^l),
\end{align*}
which shows that the sum of the coefficients of the $\beta_{0,n}$
terms is zero.  We only needed the fourth relation on the base to
establish this result.

The coefficients of the terms involving $\alpha_{0,n}$ are
\begin{align*}
[\alpha_{0,k-1},\phi_l]y^{k-1}\theta^l=
    &\alpha_{0,k+l-2}(k-1)y^{k-1}\theta^l\\
[\alpha_{0,k-1},\tau_l]y^{k-1}\zeta^l=
    &\alpha_{0,k+l-2}y^{k-1}\zeta^l\\
[\beta_{0,k-1},\gamma_{0,l}]z^{k-1}x^l=
    &\alpha_{0,k+l-2}\left(-\tfrac l2z^{k-1}x^l\right)
\end{align*}
We will show that the sum of the coefficients here satisfies
\begin{equation*}
\sum_{k+l=n+2}y^{k-1}((k-1)\theta^l+\zeta^l)-
    \tfrac l2z^{k-1}x^l=\sum_{k+l=n+3}(1-l)r_4^k\zeta^l+(k-1)r_5^k\eta^l,
\end{equation*}
and therefore it vanishes.  Thus the vanishing of the sum of these
coefficients uses only the fourth and fifth relations on the base.

\begin{align*}
0=&\sum_{k+l=n+3}(1-l)r_4^k\zeta^l+(k-1)r_5^k\eta^l\\=&
\sum_{k+l=n+3}(1-l)(s^1z^{k-1}+\sum_{i+j=k+1}\eta^i(\zeta^j+
    (i-j-1)\theta^j))\zeta^l\\
&+\sum_{k+l=n+3}(k-1)(\tfrac12 s^1y^{k-1}
    +\sum_{i+j=k+1}(1-j)\theta^i\zeta^j+\tfrac i2t^iz^{j-1})\eta^l\\
=&\!\!\!\sum_{i+j+l=n+4}\!\!(1-l)(\eta^i\zeta^j\zeta^l+
    (i-j-1)\eta^i\theta^j\zeta^l)+(i+j-2)(1-j)\theta^i\zeta^j\eta^l\\
&+s^1(\sum_{i+j+l=n+5}-(1-l)i\eta^i\eta^j\zeta^l+
    \tfrac{i+j-3}2(j-1)\eta^i\zeta^j\eta^l))\\
=&\sum_{i+j+l=n+4}(j-1)(j-1)(i+j-3)\eta^i\zeta^j\theta^l
    +j\eta^i\zeta^j\zeta^l\\
&+\sum_{i+j+l=n+5}\tfrac{(1-l)i}2s^1\eta^i\eta^j\zeta^l\\
=&\sum_{i+j+l=n+4}(j-1)(j-1)(i+j-3)\eta^i\zeta^j\theta^l+j\eta^i\zeta^j\zeta^l-
    \tfrac{il}2s^1\eta^i\eta^j\zeta^{l+1}\\
=&\sum{i+j+l=n+4}(j-1)\eta^i\zeta^j((i+j-3)\theta^l+
    \zeta^l)-\tfrac{li}2\eta^i\eta^jx^l\\
=&\sum_{k+l=n+2}y^{k-1}((k-1)\theta^l+\zeta^l)-\tfrac l2z^{k-1}x^l
\end{align*}

The coefficients of the terms involving $\gamma_{0,n}$ are
\begin{align*}
[\psi_k,\alpha_{0,l-1}]t^ky^{l-1}=&2\gamma_{k+l-1}t^ky^{l-1}\\
[\phi_k\gamma_{0,l}]\theta^kx^l=&-\gamma_{k+l-1}l\theta^kx^l\\
[\tau_k,\gamma_{0,l}]\zeta^kx^l=&\gamma_{k+l-1}\zeta^kx^l
\end{align*}
We claim that
\begin{equation*}
\sum_{k+l=n+1}2t^ky^{l-1}-l\theta^kx^l+\zeta^kx^l=
    -r_1\zeta^{n+1}+s^1r_5^{n+1}+\sum_{k+l=n+2}kr_2^k\eta^l-lr_4^kt^l
\end{equation*}
This follows from
\begin{align*}
0=&-r_1\zeta^{n+1}+s^1r_5^{n+1}+\sum_{k+l=n+2}kr_2^k\eta^l-lr_4^kt^l
\\=&
-r_1\zeta^{n+1}+\tfrac12(s^1)^2y^n
\\&+
\sum_{k+l=n+2}s^1(1-l)\theta^k\zeta^l+\tfrac k2s^1t^kz^{l-1}+\tfrac k2 s^1x^k\eta^l-ls^1z^{k-1}t^l
\\&+
\sum_{i+j+l=n+3}(i+j-1)it^i\theta^j\eta^l-2(i+j-1)t^i\zeta^j\eta^l\\&-
\sum_{i+j+l=n+3}l\eta^i\zeta^jt^l-(i-j-1)l\eta^i\theta^jt^l
\\=&
-r_1\zeta^{n+1}+\sum_{k+l=n+2}(1-l)s^1\theta^k\zeta^l
\\&+\sum_{k+i+j=n+3}(2(k-l)+i)t^i\eta^j\zeta^k-i(i+j-2)t^i\eta^j\theta^k
\\=&
s^1\zeta^1\zeta^{n+1}-t^1\eta^1\zeta^{n+1}-\sum_{k+l=n+1}ls^1\theta^k\zeta^{l+1}\\&+
\sum_{k+i+j=n+3}2(j-1)t^k\eta^i\zeta^j+(i+j-2)i\theta^kt^i\eta^j-i\zeta^kt^i\eta^j
\\=&
\sum_{k+l=n+1}2t^ky^{l-1}-l\theta^kx^l+\zeta^kx^l
\end{align*}

The demonstration of the vanishing of the coefficients of the terms
appearing in the bracket $[d^\infty,d^\infty]$ is not straightforward.
Nevertheless, this demonstration is unnecessary due to the
construction of the miniversal deformation which was given in
\cite{fp1}. We included the explicit calculations here as an
illustration of the complexity which arises in establishing the
vanishing of these coefficients by a direct calculation.

Now let us address the question of when an infinitesimal deformation
$d^i=d+\psi u +\ph \theta$, where $\psi$ is an odd and $\ph$ an even
cocycle, and $u$ is an even and $\theta$ an odd parameter, extends to
a formal deformation. Without loss of generality, we can assume that
$\psi$ is in the span of the cocycles $\xi$ and $\psi_k$, and that
$\ph$ is in the span of $\ph_k$, $\sigma_k$ and $\tau_k$, since they
differ from elements of this form by coboundaries. Suppose we write
\begin{equation*}
\psi=a_1\xi +b_1^k\psi_k,\qquad \ph=c_1^k\phi_k+g_1^k\sigma_k+h_1^k\tau_k
\end{equation*}
A morphism $f:\k[[s^1,t^k,\theta^k,\eta^k,\zeta^k]]\ra \k[[u,\theta]]$, given by
\begin{gather*}
f(s^1)=a_nu^n,\qquad
f(t^k)=b^k_nu^n\\
f(\theta^k)=c^k_nu^n\theta,\qquad
f(\eta^k)=g^k_nu^n\theta,\qquad
f(\zeta^k)=h^k_nu^n\theta
\end{gather*}
descends to one from the base $\A=\k[[s^1,t^k,\theta^k,\eta^k,\zeta^k]]/(r_m^n)$ to
$\k[[u,\theta]]$
precisely when it vanishes on the relations.  Examining the relations
carefully, we observe that the third, fourth and fifth ones only have
terms involving the product of two odd terms, and since these products
are automatically zero in $\k[[u,\theta]]$, there is nothing to check
for these relations. Substituting in the first two relations gives the
conditions
\begin{align}\label{conditions}
0=&\sum_{k+l=n}-a_kh_l^1+b_k^1g_l^1\notag\\\\
0=&\!\!\!\sum_{i+j+r=m}\!\!\!\tfrac12a_ia_jh_r^{n+1}-\!\!\!\sum_{k+l=n+2}\!\!\!ka_ib_j^kg_r^l+
\!\!\!
\sum_{\substack{i+j=m\\k+l=n+1}}\!\!\!b_i^k(kc_j^l-2h_j^l)\notag,
\end{align}
which must be satisfied for all $m$ and $n$.
It is much easier to check when the map $f$ is of degree 1, in which case the second condition breaks up into
the two separate conditions
\begin{align*}
0=&
\tfrac12(a_1)^2h_1^{n+1}-\sum_{k+l=n+2}ka_1b_1^kg_1^l\\
0=&
\sum_{k+l=n+1}b_1^k(kc_1^l-2h_1^l)
\end{align*}
The first of these two conditions coming from the second relation is
cubic and the second quadratic in the parameters $t$ and $\theta$,
while the condition derived from the first relation is also
quadratic.  Some infinitesimal deformations will not extend to second
order, because their coefficients fail the quadratic constraints,
while some will extend to second order deformations, but not to third
order, because their coefficients fail the cubic constraints.  It is
also easy to construct examples of infinitesimal deformations which
extend to a formal deformation.

For example,
$
d^i=d+\xi u +\tau_\theta
$
fails to extend to a second order deformation, because
$\tfrac12[d^i,d^i]=-\xi u\theta$. On the other hand, $d^i=d+(\xi
+\psi_1)u +\sigma_2 \theta$  extends to the second order deformation
$d_e^i=d^i-\gamma_{0,1}u\theta$, but this second order deformation
fails to extend because
$\tfrac12[d^i_e,d^i_e]=-\tfrac12\psi_1u^2\theta$.  Moreover, in the
first case, the (first) obstruction to the extension is given by
$-\xi$,  while in the second case, the (second) obstruction to the
extension is given by $-\tfrac12\psi_i$.  In fact,  it is easy to
determine the obstruction to an extension by simply plugging the
coefficients into the relations.  For example,  if
\begin{equation*}
d^i=d+\psi_1u+(\sigma_1+\phi_1+\tau_2)\theta,
\end{equation*}
then $\tfrac12[d^i,d^i]=(\xi+\psi_1-2\psi_2)u\theta$, so the
obstruction is $\xi+\psi_1-2\psi_2$. The coefficients of the cocycles
in the first obstruction are given by plugging the coefficients of the
infinitesimal extension in the quadratic parts of the relations.
However, it should be pointed out that in general, the second and
higher obstructions are not uniquely defined, because they depend on
the choice of the cochains added at each order.

For example, if $d^i=d+\xi u$, then of course, since all the relations
vanish, $d^i$ extends to a formal deformation. (In fact, $d^i$ is
itself a formal deformation.) On the other hand, the choice of the
extension to a second order deformation can affect the further
extendibility.  For example, if $d^i_e=d^i+\tau_1 u\theta$, then the
deformation cannot be extended further,  but if we let
$d^i_e=d^i+\phi_1 u\theta$, then this extension is already a
formal deformation.

What we can say is that the relations determine the maximum
extendibility of our deformation.  This is the same property as one
observes with Massey powers. The vanishing of the $n$th Massey power
means that the deformation can be extended to order $n$.  The first
nonvanishing Massey power determines the maximal order to which the
deformation can be extended.

We have to be very careful in interpreting how to use the relations,
though. In the example $d^i=d+(\xi +\psi_1)u +\sigma_2 \theta$, we can
adjust the second order extension we gave before to
$d_e^i=d^i(-\gamma_{0,1}+\tfrac12\phi_1)u\theta$, and now, the bracket
$[d_e^i,d_e^i]$ vanishes. The introduction of a cohomology class later
in the deformation corresponds to higher order terms in the polynomial
expressions of the parameters.  Here, adding the term
$\tfrac12\phi_iu\theta$ is the same as choosing
$f(\theta^1)=\tfrac12u\theta$.  Thus the computation of the
extendability of an infinitesimal deformation can be cast as follows.
Given the choice of constants $a_1$, $b^k_1$, $c^k_1$, $g^k_1$ and
$h^k_1$, do there exist constants $a_m$, $b^k_m$, $c^k_m$, $g^k_m$,
and $h^k_m$ so that the relations are satisfied? If so, then the the
deformation extends to a formal one. Thus the conditions in
\refeq{conditions} need to be solved recursively for the constants.
Moreover, the least $m$ for which a solution fails to exist determines
the maximum extendibility of the infinitesimal deformation.

Consider the formal deformation $d^f=d+\psi_1u$. Since
$\ad\tau_2(\psi_k)=2\psi_{k+1}$, we have
$(\ad\tau_2)^k(\psi_1)=2^k\psi_{k+1}$. Therefore
\begin{equation*}
d'=\exp(\ad\tau_2u)(d^f)=d+\psi_1u+\sum_{k=1}^\infty \tfrac{2^k}{k!}\psi_{k+1}u^{k+1}.
\end{equation*}
Both of these formal deformations are expressed as sums of cohomology
clases, so both of them appear as $f_*(d^\infty)$ for obvious
morphisms from the base of the miniversal deformation to
$\k[[u,\theta]]$. This example illustrates the nonuniqueness of the
morphism from the base of the versal deformation to the base
$\k[[u,\theta]]$ such that $f_*(d^\infty)\sim d^f$.
\subsection{Miniversal Deformations of Type 010}
Let $D(\ph)=\br{\ph}{\psd{0,1,1}1}$. Then we have the following table
of coboundaries.
\begin{align*}
D(\phd{1, q, n-q-1}1) &=\psd{0, 1+q, n-q}1\quad\quad D(\psd{0, p, n-p}1) =0\\
D(\phd{0, p, n-p}2) &=-\psd{0,p, n-p+1}1\quad D(\psd{1, q, n-q-1}2) =\phd{1,q,n-q}1+\phd{0,q+1,n-q}2\\
D(\phd{0, p, n-p}3) &=-\psd{0,p+1, n-p}1\quad D(\psd{1, q, n-q-1}3) =\phd{1,q+1,n-q-1}1+\phd{0,q+1,n-q}3
\end{align*}
The cohomology is given by
\begin{align*}
&H^1=\langle \psd{0,0,1}1,\psd{0,1,0}1,\phd{1,0,0}1+\phd{0,0,1}3,\phd{1,0,0}1+\phd{0,1,0}2\rangle\\
&H^n=\langle\phd{1,0,n-1}1+\phd{0,0,n}3,\phd{1,n-1,0}1+\phd{0,n,0}2\rangle \quad n>1.
\end{align*}
Let us label the cohomology classes as follows.
\begin{gather*}
\psi_1=\psd{0,0,1}1\qquad\psi_2=\psd{0,1,0}1\\
\phi_n=\phd{1,0,n}1+\phd{0,0,n+1}3\qquad \sigma_n=\phd{1,n,0}1+\phd{0,n+1,0}2\qquad n\ge 0
\end{gather*}
The universal infinitesimal deformation is given by
\begin{equation*}
d^1=\psd{0,1,1}1+\psi_1 t^1 +\psi_2t^2+\phi_n \theta^n +\sigma_n\eta^n,
\end{equation*}
where $t^i$ are even parameters and $\theta^n$ and $\eta^n$ are odd
parameters. Let
\begin{align*}
\alpha_{k,l}=\phd{0,k,l}3\qquad \beta_{k}=\phd{0,0,k}2\qquad
\tau_{k,l}=\psd{1,k,l}2\qquad\xi_{k,l}=\psd{1,k,l}3
\end{align*}
These cochains are preimages of a basis of the coboundaries, so it is
possible to express the miniversal deformation in the form
\begin{equation*}
d^\infty=d^1+\alpha_{k,l}x^{k,l}+\beta_ky^k+\tau_{k,l}u^{k,l}+\xi_{k,l}v^{k,l}.
\end{equation*}
It turns out that we do not need all of the above cochains to
construct the miniversal deformation. Let us denote
\begin{equation*}
\gamma_k=\alpha_{k,1}\qquad\alpha_{0,k}=\epsilon_k\qquad k>0
\end{equation*}
and set $r^k=x^{k,1}$, $s^k=x^{0,l}$. Then we will show that the miniversal deformation
can be expressed in the form
\begin{equation*}
d^\infty=d^1+\gamma_kr^k+\epsilon_ks^k+\beta_ky^k+\tau_{k,l}u^{k,l}+\xi_{k,l}v^{k,l},
\end{equation*}
where $\gamma_k$, $\epsilon_k$ and $\beta_k$ are even cochains defined
for $k\ge 1$, $\tau_{k,l}$ and $\xi_{k,l}$ are odd cochains defined
for $k,l\ge 0$, and actually, $v^{k,0}=0$.

The brackets we need to compute in order to determine $[d^1,d^1]$ are
\begin{align*}
[\psi_1,\phi_k]=&[\psi_2,\sigma_k]=[\phi_k,\sigma_l]=0\\
[\psi_1,\sigma_k]=&
\begin{cases}
-\psi_1&k=0\\
D(\epsilon_1)&k=1\\
D(\gamma_{k-1})&k>1
\end{cases}\\
[\psi_2,\phi_k]=&
\begin{cases}
-\psi_2&k=0\\
D(\epsilon_k)& k>0
\end{cases}\\
[\phi_k,\phi_l]=&\phi_{k+l}(k-l),\quad
[\sigma_k,\sigma_l]=\sigma_{k+l}(k-l)
\end{align*}
The second order relations are
\begin{gather*}
t^1\eta^0=t^2\theta^0=
\tfrac12\!\!\sum_
{\substack{
k+l=n
}}\!\!
(k-l)\theta^k\theta^l=
\tfrac12\!\!\sum_
{\substack{
k+l=n
}}\!\!
(k-l)\eta^k\eta^l=0
\end{gather*}
The second order deformation is given by
\begin{equation*}
d^2=d^1+\epsilon_1t^1\eta^1+\epsilon_{k}t^2\theta^k +\gamma_kt^1\eta^{k+1}.
\end{equation*}
We next compute the brackets which are necessary to compute
$[d^2,d^2]$. These are the brackets of cohomology classes and the
$\gamma$ and $\epsilon$ terms. Note that since these terms first
appeared in the brackets of cohomology classes, their coefficients in
the miniversal deformation have order two, so the brackets below have
order 3.
\begin{align*}
[\psi_1,\gamma_k]=
&\begin{cases}
-D(\epsilon_1)&k=1\\
-D(\gamma_{k-1})& k>1\\
\end{cases}&
[\psi_1,\epsilon_k]=&
\begin{cases}
\psi_1&k=1\\
-D(\beta_{k-1})&k>1
\end{cases}\\
[\psi_2,\gamma_k]=&0&[\psi_2,\epsilon_k]=&0\\
[\phi_k,\gamma_l]=&D(\xi_{l-1,k})k&[\phi_k,\epsilon_l]=&\phi_{k+l-1}k+\epsilon_l(k-l)\\
[\sigma_k,\gamma_l]=&-\gamma_{k+l}l&[\sigma_k,\epsilon_l]=&0
\end{align*}
Let us show that the terms involving $\gamma$ and $\epsilon$ cochains
cancel, at least up to fourth order.
Consider the following terms
\begin{align*}
[\sigma_{k+1},\gamma_l]\eta^{k+1}t^1\eta^{l+1}=&-l\gamma_{n+1}t^1\eta^{k+1}\eta^{l+1},k+l=n\\
[\gamma_k,\sigma_{l+1}]t^1\eta^{k+1}\eta^{l+1}=&k\gamma_{n+1}t^1\eta^{k+1}\eta^{l+1},k+l=n\\
[\sigma_0,\gamma_{n+1}]t^1\eta^0\eta^{n+2}=&-(n+1)\gamma_{n+1}t^1\eta^0\eta^{n+2}.
\end{align*}
Summing the first two types and dividing by $\tfrac12$ and adding the
last term gives the second order relation involving $\eta$ cochains,
plus the term $t^1\eta^0\eta^{n+2}$. But this term is zero up to
fourth order, using the second order relation $t^1\eta^0=0$. The terms
involving $\epsilon$ cochains are handled similarly. The third order
relations are
\begin{gather*}
t^1\eta^0-t^1(t^1\eta^1+t^2\theta^1)=t^2\theta^0=0\\
nt^1\theta^n\eta^1+\!\!\sum_
{\substack{
k+l=n
}}\!\!\tfrac12
(k-l)\theta^k\theta^l+kt^1\theta^k\theta^{l+1}=
\tfrac12\!\!\sum_
{\substack{
k+l=n
}}\!\!
(k-l)\eta^k\eta^l=0
\end{gather*}
Notice that only two of them have been modified from the second order
relations. The third order deformation is
\begin{equation*}
d^3=d^2-\epsilon_1(t^1)^2\eta^2-\gamma_k(t^1)^2\eta^{k+2}-\beta_kt^1t^2\theta^{k+1}-\xi_{k,l}lt^1\theta^l\eta^{k+2}.
\end{equation*}
Next,  we compute all brackets whose order is 4. Note that the
brackets of $\gamma$ and $\epsilon$ terms appear here, because they
have order 4, and so play no role in the construction of the third
order deformation.
\begin{align*}
[\gamma_k,\gamma_l]=&0\\
[\gamma_k,\epsilon_l]=&\alpha_{k,l}(1-l)\\
[\epsilon_k,\epsilon_l]=&\epsilon_{k+l-1}(k-l)\\
[\psi_1,\beta_k]=&0\\
[\psi_2,\beta_k]=&\begin{cases}
\psi_1&k=1\\
-D(\beta_{k-1})&\text{otherwise}
\end{cases}\\
[\phi_k,\beta_l]=&-\beta_{k+l}l\\
[\sigma_k,\beta_l]
=&
\begin{cases}
\beta_l&k=0\\
2D(\tau_{0,l-1})-\phi_l+\epsilon_{l+1}&k=1\\
D(\tau_{k-1,l-1})(k+1)-D(\xi_{k-2,l})+\alpha_{k-1,l+1}&k>1
\end{cases}\\
[\psi_1,\xi_{k,l}]=&
\begin{cases}
\phi_l&k=0\\
D(\xi_{k-1,l})&\text{otherwise}
\end{cases}\\
[\psi_2,\xi_{k,l}]=&
\begin{cases}
\gamma_{k+1}&l=1\\
\alpha_{k+1,l}&l>1
\end{cases}\\
[\phi_n,\xi_{k,l}]=&\xi_{k,n+l}(n-l)\\
[\sigma_n,\xi_{k,l}]=&-\xi_{k+n,l}(k+1)
\end{align*}
Note the appearance of the terms $\alpha_{n,k}$. Let us show that the
coefficient of such terms is zero, up to order 5. We have
\begin{align*}
[\gamma_k,\epsilon_l]t^1t^2\eta^{k+1}\theta^l=&\alpha_{k,l}(1-l)t^1t^2\eta^{k+1}\theta^l\\
[\sigma_{k},\beta_l](-t^1t^2\eta^{k}\theta^{l+1})=&\alpha_{k-1,l+1}(-t^1t^2\eta^{k}\theta^{l+1})\\
[\psi_2,\xi_{k,l}](-lt^1t^2\theta^l\eta^{k+2})=&\alpha_{k+1,l}(-lt^1t^2\theta^l\eta^{k+2})
\end{align*}
Adjusting the indices, and interchanging the odd terms on the third
equation, one sees that these terms add up to zero.

In the fourth order deformation, it is necessary to introduce the
cochains $\tau_{k,l}$, which therefore will be of order 4. Some
modifications to the relations occur, and some additional $\gamma$,
$\beta$ and $\xi$ terms will be added.  We will not give the fourth
order deformation explicitly here, because we will compute the
miniversal deformation directly by recursion. In order to do so,  we
need to compute all the brackets of all of the remaining terms with
each other.
\begin{align*}
[\gamma_k,\beta_l]=&
\begin{cases}
-D(\tau_{0,l-1})l+\phi_l+\epsilon_l(1-l)&k=1\\
-D(\tau_{k-1,l-1})l+D(\xi_{k-2,l})l
-\alpha_{k-1,l+1}(k-l)
&k>1
\end{cases}\\
[\epsilon_k,\beta_l]=&-\beta_{k+l-1}l\\
[\xi_{k,l},\gamma_n]=&-\xi_{k+n,l}(1-l)\\
[\xi_{k,l},\epsilon_n]=&-\xi_{k,l+n}(n-l)\\
[\xi_{k,l},\beta_n]=&-\tau_{k,l+n-1}+\xi_{k-1,l+n}k\\
[\psi_1,\tau_{k,l}]=&
\begin{cases}
\beta_{l+1}&k=0\\
D(\tau_{0,l})-\phi_{l+1}+\epsilon_{l+2}&k=1\\
D(\tau_{k-1,l})-D(\xi_{k-2,l+1})+\alpha_{k-1,l+2}&\text{otherwise}
\end{cases}\\
[\psi_2,\tau_{k,l}]=&
\begin{cases}
\sigma_k&l=0\\
D(\tau_{k,l-1})&\text{otherwise}
\end{cases}\\
[\tau_{k,l},\phi_n]=&\tau_{k,n+l}(l+1)\\
[\tau_{k,l},\sigma_n]=&-\tau_{k+n,l}(n-k)\\
[\tau_{k,l},\gamma_n]=&\tau_{k+n,l}-\xi_{k+n-1,l+1}n\\
[\tau_{k,l},\epsilon_n]=&\tau_{m,l+n-1}\\
[\tau_{k,l},\beta_n]=&\tau_{k-1,l+n}k\\
[\beta_k,\beta_l]=&[\tau_{k,l},\tau_{m,n}]=[\tau_{k,l},\xi_{m,n}]=[\xi_{k,l},\xi_{m,n}]=0
\end{align*}
Let us collect the terms involving the coboundaries of the $\gamma$ cochains.
Including the coefficients, we obtain
\begin{align*}
[\psi_1,\sigma_k]t^1\eta^k=&D(\gamma_{k-1})t^1\eta^k\\
[\psi_1,\gamma_k]t^1r^k=&-D(\gamma_{k-1})t^1r^k\\
[d,\gamma_k]r^k=&-D(\gamma_k)r^k
\end{align*}
Since the sum of all terms involving the same index in $\gamma$ must
vanish, we obtain the recursive relation
$
r^k=t^1\eta^{k+1}-t^1r^{k+1}
$,
 from which it follows that
\begin{equation*}
r^k=t^1\eta^{k+1}+\sum_{n=1}^\infty\s n (t^1)^{n+1}\eta^{k+n+1}.
\end{equation*}

For the $\epsilon$ cochains we have
\begin{align*}
[\psi_1,\sigma_1]t^1\eta^1=&D(\epsilon_1)t^1\eta^1\\
[\psi_2,\phi_k]t^2\theta^k=&-D(\epsilon_k)t^2\theta^k\\
[d,\epsilon_k]s^k=&-D(\epsilon_k)s^k.
\end{align*}
which yields
\begin{align*}
s^1=&t^2\theta^1+t^1\eta^1\\
s^k=&t^2\theta^k, k>1
\end{align*}

The terms involving coboundaries of $\beta$ cochains are
\begin{align*}
[\psi_1,\epsilon_k]t^1s^kx=&-D(\beta_{k-1})t^1s^k\\
[\psi_2,\beta_k]t^2y^k=&-D(\beta_{k-1})t^2y^k\\
[d,\beta_k]y^k=&-D(\beta_k)y^k.
\end{align*}
Since we have $y^k=-t^1s^{k+1}-t^2y^{k+1}$, it follows easily that
\begin{equation*}
y^k=-t^1t^2\sum_{n=0}^\infty\s n(t^2)^n\theta^{n+k+1}.
\end{equation*}
The terms involving $\tau_{k,l}$ are more complicated.
\begin{align*}
[\sigma_k,\beta_l]\eta^ky^l=&D(\tau_{k-1,l-1})(k+1)\eta^ky^l\\
[\psi_1,\tau_{k,l}]t^1u^{k,l}=&D(\tau_{k-1,l})t^1u^{k,l}\\
[\psi_2,\tau_{k,l}]t^2u^{k,l}=&D(\tau_{k,l-1})t^2u^{k,l}\\
[\gamma_k,\beta_n]r^ky^n=&-D(\tau_{k-1,n-1})nr^ky^n\\
[d,\tau_{k,l}]u^{k,l}=&-D(\tau_{k,l})u^{k,l}.
\end{align*}
This yields the following recursion relation
\begin{equation*}
u^{k,l}=((k+2)\eta^{k+1}-(l+1)r^{k+1})y^{l+1}+t^1u^{k+1,l}+t^2u^{k,l+1},
\end{equation*}
which gives a power series expression for $x^{k,l}$.
Finally, the terms involving $\xi$'s are
\begin{align*}
[\ph_n,\gamma_k]\theta^nr^k=&D(\xi_{k-1,n})n\theta^nr^k\\
[\sigma_k,\beta_l]\eta^ky^l=&-D(\xi_{k-2,l})\eta^ky^l\\
[\psi_1,\tau_{k,l}]t^1u^{k,l}=&-D(\xi_{k-2,l+1})t^1u^{k,l}\\
[\psi_1,\xi_{k,l}]t^1v^{k,l}=&D(\xi_{k-1,l})t^1v^{k,l}\\
[\gamma_k,\beta_n]r^ky^n=&D(\xi_{k-2,n})nr^ky^n\\
[d,\xi_{k,l}]v^{k,l}=&-D(\xi_{k,l})v^{k,l},
\end{align*}
from which we deduce that
\begin{equation*}
v^{k,l}=l\theta^lr^{k+1}-\eta^{k+2}y^l+lr^{k+2}y^l-t^1u^{k+2,l-1}+t^1v^{k+1,l}.
\end{equation*}
Now let us study the relations on the base.  These are the
coefficients of the cocycles. There are three terms involving
$\psi_1$.  From $[\psi_1,\epsilon_1]$ we obtain $t^1s^1$, from
$[\psi_2,\beta_1]$ we obtain $t^2y^1$, and from $[\psi_1,\sigma_0]$ we
obtain $-t^1\theta^0$. Thus the coefficient of $\psi_1$ is
$
-t^1\theta^0+t^2y^1+t^1s^1
$.
The only term involving $\psi_2$ is $[\psi_2,\theta_0]$, giving $t^2\theta^0=0$.

For $\phi_n$ we obtain the terms
$\tfrac12[\phi_k,\phi_l]=\phi_n\tfrac12(k-l)\theta^k\theta^l$, where
$k+l=n$. From the brackets $[\phi_k,\epsilon_l]$, if we re-index in
the form $[\phi_k,\epsilon_{l+1}]$, so that we can sum for all $k$,
$l$, not just $l>0$, and then sum the corresponding terms with $k$ and
$l$ interchanged, we obtain $\tfrac12(k-l)\theta^ks^{l+1}\phi_n$, for
$k+l=n$. Similarly, from  $[\gamma_1,\beta_{n}]$, we get
$-\tfrac12(k-l)r^1y^nn$. Lastly, there are the terms
$[\sigma_1,\beta_n]$, contributing $-\eta^1y^n\phi_n$,
$[\psi_1,\tau_{1,n-1}]$, giving $-t^1u^{1,n-1}\phi_n$, and
$[\psi_1,\xi_{0,n}]$, yielding $t^1v^{0,n}\phi_n$. These last terms
are only defined when $n\ge1$.  Putting this altogether, we obtain
\begin{equation*}
-\eta^1y^n-t^1u^{1,n-1}+t^1v^{0,n}-s^1y^{n}+\tfrac12\!\!\sum_{k+l=n}\!\!(k-l)(\theta^k\theta^l+\theta^ks^{l+1})
\end{equation*}

Finally, let us examine the terms involving the $\sigma$ cochains.
From $\tfrac12[\sigma_k,\sigma_l]$, we obtain
$\sigma_n\tfrac12(k-l)\eta^k\eta^l$, when $k+l=n$. From
$[\psi_2,\tau_{n,0}]$, we obtain $t^2u^{n,0}$. Thus the corresponding
relation is $$t^2u^{n,0}+\tfrac12\sum_{k+l=n}(k-l)\eta^k\eta^l.$$

Putting these all together, the relations on the base of the
miniversal deformation are
\begin{gather*}
t^2\theta^0=0\\
-t^1\theta^0+t^1s^1+t^2y^1=0\\
-\eta^1y^n-t^1u^{1,n-1}+t^1v^{0,n}-s^1y^{n}+\tfrac12\sum_{k+l=n}(k-l)(\theta^k\theta^l+\theta^ks^{l+1})=0\\
t^2u^{n,0}+\tfrac12\sum_{k+l=n}(k-l)\eta^k\eta^l=0
\end{gather*}
Note that the $s$, $y$, $u$ and $v$ coefficients above can be
expressed in terms of the parameters $t^1$, $t^2$, $\theta$ and
$\eta$. The first two relations are odd, and the last two are even, so
only the first two play a role in determining whether an infinitesimal
deformation extends to a formal one. Substituting for the $s^1$ and
$y^1$ coefficients in the second relation gives
\begin{equation*}
-t^1(\theta^0+t^2\theta^1+(t^1)^2\eta^1-\sum_{n=0}^\infty(t^2)^{n+2}\theta^{n+2})=0,
\end{equation*}
so that $t^1=0$ certainly solves this equation. Now suppose that
$d^i=d+\psi u +\ph\theta$ is an infinitesimal deformation of $d$. Let
\begin{equation*}
t^i=a^i_ku^k,\qquad
\theta^i=b_k^iu^k\theta,\qquad \eta^i=c_k^iu^k\theta,
\end{equation*}
where $\psi=a^1_1\psi_1+a^1_1\psi_2$ and $\ph=b^i_1\phi_i
+c^i_1\sigma_i$. Our goal is to solve for coefficients $a^i_k$,
$b^i_k$, and $c^i_k$ so that the first two relations are satisfied. If
$a^1_1=a^2_1=0$, then by choosing $a^i_k=0$ for all $k$, we obtain a
solution. The relations transform to
\begin{gather*}
\sum_{i+j=m}a^2_ib^0_j=0\\
-\sum_{i+j=m}a^1_ib^0_j+\sum_{i+j+k=m}\!\!\!\!\!\!\!\!a^1_i(a^2_jb^1_k+a^1_jc^1_k)+
\!\!\!\!\!\!\!\!\!\sum_{\substack{n=0\cdots\infty\\i+j+k_1\mplus k_{n+2}=m}}
\!\!\!\!\!\!\!\!
\s n a_i^1c_j^{n+2}\prod_{l=0}^{n+2}a^2_{k_l}=0
\end{gather*}
Suppose that $a_k^1\ne 0$ for some $k$. Then the second relation
transforms into the simpler
\begin{equation*}
-b^0_m+\sum_{j+k=m}\!\!\!\!\!(a^2_jb^1_k+a^1_jc^1_k)+
\!\!\!\!\!\!\!\!\!\sum_{\substack{n=0\cdots\infty\\j+k_1\mplus k_{n+2}=m}}\!\!\!\!\!\!\!\!
(-1)^nc_j^{n+2}\prod_{l=0}^{n+2}a^2_{k_l}=0.
\end{equation*}
Otherwise, the second relation is satisfied automatically. Looking at
this transformed second relation, we observe that $b_1^0=0$, since
this is the only term in the expression for $m=1$.  If $t^2=0$, then
the first relation is satisfied, and we can always solve for
coefficients $b_m^0$ to satisfy the second relation, for arbitrary
choices of the other coefficients. On the opposite end of the
spectrum, if $a^2_1\ne 0$, then $b_k^0=0$ for all $k$, or we encounter
the obstruction $a^2_1b^0_m\psi_1$ at the first level $m$ for which
$b^0_m\ne 0$. We may encounter an earlier obstruction coming from the
second relation, because, for example, for $m=2$, we obtain from the
second equation the relation $b_2^0+a^2_1b^1_1+a^1_1c^1_1=0$.
Actually, in the case that $a^1_1\ne 0$, this is the only obstruction
coming from the second relation, because the term $a_1^1c_m^1$ appears
in the sum for degree $m+1$, which means that for a certain choice of
the coefficient $c_m^1$, one can make the sum equal to zero. Thus, for
an infinitesimal deformation, we can resolve the extendibility to a
formal deformation as follows.

If $a^1_1=a^2_1=0$, then the deformation extends trivially. In fact,
the infinitesimal deformation is already a formal deformation.
Otherwise, if $b_1^0\ne 0$, then the deformation does not extend to
second order, and the obstruction is
$a^1_1b_1^0\psi_1+-a^2_1b_1^0\psi_2$. If $b_1^0=0$, then the
deformation does extend to second order. If $a_1^1=0$, then the second
relation vanishes by choosing $a_k^1=0$ for all $k$, and the first
relation vanishes by choosing $b_k^0=0$ for all $k$, so the
deformation extends to a formal one. On the other hand, if $a_1^2=0$,
then by choosing $a_k^2=0$ for all $k$, we can make the first relation
vanish, and by choosing $b_m^0$ appropriately, we can make the second
relation vanish for each order $m$. Thus the deformation extends to a
formal one.  When neither $a_1^1$ nor $a_1^2$ vanish, we must have
$b_k^0=0$ for all $k$. If $a_1^2b_1^1+a_1^1c_1^1\ne 0$, then no second
order deformation extends to third order. On the other hand, if it
does vanish, then there is always a choice of the coefficients $c_k^1$
such that the deformation extends to a formal one.

\section{Conclusions}
The notion of miniversal deformations has been around for quite some
time.  Nevertheless, there were some confusions in the early
literature, as was mentioned in \cite{ff2}. We have felt that some of
the confusion arises because of the lack of concrete examples. Our
purpose in this article has been to give some explicit constructions
of miniversal deformations, and use them to address  the classical
question of when an $n$-th order deformation extends to a formal one.

In \cite{ff2}, relations on the base of the miniversal deformation
were found, without actually constructing a miniversal deformation.
Since these relations alone determine the extendibility of a
deformation, it is interesting to note that they can sometimes be
found without having to carry out the complete construction.

Another important classical question is: 

\emph{Given that an extension to a
formal deformation exists, how can you determine the equivalence
classes of nonequivalent extensions?}

We did not address this problem
in this article. In some of our other work \cite{fp4,fp5}, we have
addressed the problem of how to extend a codifferential of degree $n$
to a more general \linf\ structure, with possibly infinitely many
terms. We studied the equivalence classes of these extensions, and the
classification problem is quite tricky. Since these extensions can be
thought of as specializations of deformations of the \linf\ structure
determined by the degree $n$ codifferential, where the even parameters
are given fixed values, and the odd parameters have been set equal to
zero, there is a close relationship between the problem of
classification of extensions and the classification of deformations up
to equivalence. Thus, we don't expect the classification of
deformations to be easy.


\providecommand{\bysame}{\leavevmode\hbox to3em{\hrulefill}\thinspace}
\providecommand{\MR}{\relax\ifhmode\unskip\space\fi MR }
\providecommand{\MRhref}[2]{%
  \href{http://www.ams.org/mathscinet-getitem?mr=#1}{#2}
}
\providecommand{\href}[2]{#2}

\end{document}